\documentclass{daj}

\usepackage{graphicx}
\usepackage{amssymb}
\usepackage{epstopdf}
\usepackage{amsmath,amsthm,amscd,amssymb}
\usepackage{latexsym}
\usepackage{hyperref}

\theoremstyle{plain}
\newtheorem{theorem}{Theorem}[section]
\newtheorem{lemma}[theorem]{Lemma}

\theoremstyle{definition}
\newtheorem{definition}[theorem]{Definition}

\newtheorem{question}[theorem]{Question}

\newtheorem{case[theorem]}{Case}

\theoremstyle{remark}
\newtheorem{remark}[theorem]{Remark}

\numberwithin{equation}{section}

\dajAUTHORdetails{%
  title = {Gabor orthogonal bases and convexity}, 
  author = {Alex Iosevich and Azita  Mayeli},
  plaintextauthor = {Alex Iosevich and Azita  Mayeli},
  plaintexttitle = {Gabor orthogonal bases and convexity}, 
  runningtitle = {Gabor orthogonal bases and convexity}, 
  runningauthor = {Alex Iosevich and Azita  Mayeli},
  copyrightauthor = {Alex Iosevich and Azita  Mayeli},
  keywords = {Gabor basis, Spectral sets, Zeroes of the Fourier transform, Non-vanishing Gaussian curvature},
}   

\dajEDITORdetails{%
   year={2018},
   number={19},
   received={7 September 2017},   
   published={12 December 2018},  
   doi={10.19086/da.5952},       
}   

\begin{document}

\begin{frontmatter}[classification=text]

\title{Gabor Orthogonal Bases and Convexity} 

\author[iosevich]{Alex Iosevich\thanks{Partially supported by the NSA Grant H98230-15-1-0319.}}
\author[mayeli]{Azita Mayeli\thanks{Partially supported  by  PSC-CUNY grant 69625-00-47 and 48, jointly funded by The Professional Staff Congress and The City University of New York.}}

\begin{abstract}
Let $g(x)=\chi_B(x)$ be the indicator function of a bounded convex set  $B$ in $\mathbb R^d$, $d\geq 2$,  with a smooth boundary and everywhere non-vanishing Gaussian curvature. Using a combinatorial approach we prove that if $d \neq 1 \mod 4$, then there does not exist $S \subset {\Bbb R}^{2d}$ such that ${ \{g(x-a)e^{2 \pi i x \cdot b} \}}_{(a,b) \in S}$ is an orthonormal basis for $L^2({\Bbb R}^d)$.
\end{abstract}
\end{frontmatter}

\section{Introduction}
 The basic question we ask in this paper is, which functions $g$ can serve as window functions for orthogonal Gabor bases for $L^2({\Bbb R}^d)$? 

Let $g\in L^2({\Bbb R}^d)$ and $S\subset \Bbb R^{2d}$ be a countable subset. The Gabor system associated with $g$ and $S$ is defined to be the set of functions
$$\mathcal G(g, S) = \{g(x-a) e^{-2 \pi i x \cdot b} \}_{(a,b) \in S}. $$

\begin{definition} We say that ${\mathcal G}(g,S)$ is an orthonormal basis for $L^2({\Bbb R}^d)$ with the window function $g$, $\|g\|=1$,   and the spectrum $S$ if ${\mathcal G}(g,S)$  is  complete  in $L^2({\Bbb R}^d)$ and the   vectors are mutually orthogonal in the sense that 
\begin{equation} \label{orthoformula} \int g(x-a) g(x-a') e^{-2 \pi i x \cdot (b-b')} dx=0 \ \text{for all} \ (a,b) \not=(a',b'). \end{equation} 
\end{definition} 

The theory of Gabor bases and frames has undergone much development in recent decades. But as Gr\"ochenig points out in his seminar article (\cite{Gr14}), "there has been little progress on the original question of how to determine which windows and lattices generate a Gabor frame". In this paper we take a small step in this direction in the context of general Gabor orthogonal bases, not necessarily lattices. The question of which functions $g$ can serve as window function for orthogonal Gabor bases is typically studied using the Balian-Low theorem and its variants. See, for example, \cite{B88}, \cite{CLMP17}, \cite{D90}, \cite{GM13}, \cite{NO17} and the references contained therein. See also \cite{AAK17}, \cite{HW04} and \cite{LW03} for closely related results. However, in this paper our aim is to rule out a class of window functions which are indicator functions of bounded sets. In such cases, the Fourier transform of the window function is instantly poorly localized and hence Balian-Low type theorems are difficult to utilize. 

When $S=A \times B$ and $g(x)=\chi_E(x)$, with $E$ a bounded subset of ${\Bbb R}^d$ with non-zero measure, then it is not difficult to see that $E$ must tile ${\Bbb R}^d$ by translation, with $A$ serving as a tiling set. Similarly, in this case ${\{e^{-2 \pi i x \cdot b} \}}_{b \in B}$ would have to be an orthogonal basis for $L^2(E)$. This largely reduces the orthogonal Gabor basis problem to the tiling and orthogonal exponential basis components. As a result, we can use the theory of tiling and orthogonal exponential bases to rule out the possibility that $\chi_E$, for a given set $E$, is a window function for an orthogonal Gabor basis. 

This allows us to rule out the possibility that $\chi_E$ is the window function for an orthogonal Gabor basis with $S=A \times B$ for many classes of sets $E$. For example, if $E=B_d$ is the unit ball, we can rule out $\chi_{B_d}$ in two different ways. First of all, the unit ball does not tile by translation. It is known that a convex body tiles by translation only if it is a polyhedron of a certain type. See, for example, \cite{M80}. Another way to prove that $\chi_{B_d}$ is not a window function in the case $S=A \times B$ is by showing that $L^2(B_d)$ does not possess an orthogonal exponential basis. This was established by the first listed author, Nets Katz and Steen Pedersen in \cite{IKP01}. 

When $S$ is not of the form $A \times B$, the problem of determining the possible window functions becomes considerably more difficult. In this paper we rule out the possibility that the window function $g$ is an indicator function of a symmetric convex set with a smooth boundary and everywhere non-vanishing Gaussian curvature. The approach used by the first listed author, Katz and Pedersen (\cite{IKP01}) or the one employed the first listed author, Katz and Tao in \cite{IKT01} are difficult to apply here owing to the fact that the Gabor spectrum is not assumed to be of the form $A \times B$. However, the geometric approach used in \cite{IR03} combined with a suitable combinatorial pigeon-holing technique allows us to rule out a large class of window functions. Our main result is the following. 

\vskip.125in 

\begin{theorem} \label{main} Let $K$ denote a bounded convex set in $\Bbb R^d$, $d \neq 1 \mod 4$, symmetric with respect to the origin. Suppose that 
$\partial K$ is smooth and has everywhere non-vanishing Gaussian curvature. Then there does not exist an orthonormal Gabor basis for $L^2({\Bbb R}^d)$ with the window function $g(x)=|K|^{-1/2}\chi_K(x)$. 


\end{theorem} 

\vskip.125in 

\begin{remark} The case $d=1\mod 4$ requires a fundamentally different approach. This issue is taken up in \cite{IKMN18}. It is interesting to note that a similar issue arose in the context of orthogonal sets of exponentials in \cite{IR03}. \end{remark} 

\begin{remark} It should not be difficult to modify the proof of Theorem \ref{main} to handle the more general case when $\partial K$ is smooth and the boundary contains at least one point where the Gaussian curvature does not vanish. We shall address these issues in the sequel. Other related issues are raised in Section \ref{openproblems} below. 
\end{remark} 

\vskip.125in 


\vskip.125in 


\section{Expansion estimates}

\section{Proof of Theorem \ref{main}} 

The structure of our argument is the following. 

\begin{itemize} 
\item i) First we use the general theory of Gabor frames to prove that $S$ is well-distributed in the sense that ${\Bbb R}^{2d}$ can be tiled by cubes of side-length $C>0$ such that every cube contains at least one pair $(a,b) \in S$. This argument is based purely on the completeness of Gabor frames.

\vskip.125in 

\item ii) We use a pigeon-holing argument to show that there exist at least three points $(a_j,b_j) \in S$, $j=1,2,3$, such that $|a_i-a_j|$ is small and $|b_i-b_j|$, $i \not=j$, is suitably large. 

\vskip.125in 

\item iii) Finally, we use orthogonality and an asymptotic formula for the Fourier transform of the indicator function of a symmetric convex body with a smooth boundary and everywhere non-vanishing Gaussian curvature to see that the triplet in item ii) cannot exist, obtaining a contradiction. This is where the assumption $d \not=1 \mod 4$ plays a key role. 

\end{itemize} 

\vskip.125in 

\subsection{Basic structure of orthonormal Gabor bases} 

\vskip.125in 

We need the following basic structure theorem. 

\begin{lemma} \label{structuretheorem} Let ${\mathcal G}(g,S)$ be an orthonormal basis for $L^2({\Bbb R}^d)$, where $g$ is in  $L^2({\Bbb R}^d)$ with norm $1$. Then the following hold: \begin{itemize} 

\vskip.125in 

\item i) There exists $c>0$ such that for any $(a,b), (a',b') \in S$, $(a,b) \not=(a',b')$, 
$$|a-a'|+|b-b'| \ge c.$$

\vskip.125in 

\item ii) There exists $C>0$ such that any cube of side-length $C$ in ${\Bbb R}^{2d}$ contains at least one point of $S$. 

\vskip.125in 

\item iii) $S$ has uniform density equal to 1: $D^{+}(S)= D^{-}(S)=1$. 
\end{itemize}
\end{lemma} 

\vskip.125in 

\subsection{Proof of part i)} Observe that if $(a,b)=(a',b')$, then the orthogonality relation takes the form 
 \begin{align}\notag\label{eq:1}
  \int {|g(x)|}^2 dx=1.\end{align}

The continuity of the integral implies  that 
 if $(a,b)$ is sufficiently close to $(a',b')$, then 
 
$$ \int g(x-a) \overline{g(x-a')} e^{-2 \pi i x \cdot (b-b')} dx \not=0.$$ 
 The claim follows.

\vskip.125in 

\subsection{Proof of part ii)} The proof follows from the following result. 
 Define the short time Fourier transform (STFT) (see e.g. \cite{G01}) by 
$$ V_gf(t,\nu)=\int f(x) \overline{g(x-t)} e^{-2 \pi i x \cdot \nu} dx, \ (t,\nu)\in \Bbb R^{2d} .$$ 

\begin{lemma} \label{besselidentity} (\cite{AAK17})
Given $g\in L^2(\Bbb R^d)$ with $\|g\|=1$ and a countable set $S\subset \Bbb R^{2d}$, the system ${\mathcal G}(g,S)$ is an orthonormal basis for $L^2(\Bbb R^d)$ if and only if  for any $f\in L^2(\Bbb R^d)$ with $\|f\|=1$, 
\begin{equation} \label{besselformula} \sum_{\alpha \in S} {|V_g f(w-\alpha)|}^2=1 \ \text{for almost all} \  w \in {\Bbb R}^{2d}.\end{equation} 
\end{lemma} 

In other words,   $|V_gf|^2+S=1$ is a tiling, hence $S$ has asymptotic density $1$ and 
the conclusion of (ii) follows from \cite{KL96}. 

\vskip.125in 

\subsection{Proof of part iii)}  The proof follows from  Corollary 4 of \cite{Ramanathan-Steger94}.  

\vskip.125in 

We shall also make use of the following result (\cite{K18}). 

\begin{theorem} \label{localizedgabor} Let $\Omega$ be a bounded open set. Suppose that ${\mathcal G}(g,S)$ be an orthonormal basis for $L^2({\Bbb R}^d)$, where $g=\chi_K$, with $K$ a bounded set of volume $1$. Then, in particular, ${\mathcal G}(g,S)$ is complete in 
$L^2(\Omega)$. Define 
$$\pi_1(S)=\{a: (a,b) \in S\}, \quad \text{and} \quad   \pi_2(S)=\{b: (a,b) \in S\}.$$ 

Let $B$ be a small ball and $D \subset \Omega$ such that $D+B \subset \Omega$. Then 
$$ \pi_2(S \cap (\Omega+B) \times {\Bbb R}^d))$$ has positive lower density. 
\end{theorem} 

\vskip.125in 

To prove Theorem \ref{localizedgabor}, observe that orthogonality implies that for any $\phi \in L^2({\Bbb R}^d)$, with $||\phi||=1$, we have
$$ \sum_{(a,b) \in S} {\left| <\phi(t), g(t-a)e^{2\pi i b t}>\right|}^2 \le 1, $$
which we rewrite as
$$ \sum_{(a,b) \in S} {\left|V_g\phi(a,b)\right|}^2 \le 1, $$
where
$$ V_gf(a,b) := <f,  g(\cdot-a)e^{2\pi i b. \cdot}>= \int_{{\Bbb R}^d} f(t) \overline{g(t-a)} e^{-2\pi i b t} dt $$

is the short time Fourier transform (STFT for short) with $g$ as the window function. 
Apply this now to any time-frequency translate of $\phi$, i.e.   $\phi(\cdot - x)e^{2\pi i y\cdot}$,  to get
\begin{equation} \label{packing} \sum_{(a,b) \in S} {\left|V_g\phi((a,b)-(x,y))\right|}^2 \le 1, \end{equation}
which is  valid for any $\phi \in L^2({\Bbb R}^d)$. In other words, the inequality (\ref{packing}) means that  the $S$-time-frequency translates of the function ${|V_g\phi|}^2$ are {\em packing}  for $\Bbb R^{2d}$ at level one.

By our assumption of completeness in $L^2(\Omega)$, whenever the $(x,y)$-time-frequency translate of $\phi$ is supported in $\Omega$ , then the inequality  \eqref{packing} becomes an equality.

Take $\phi$ to be supported on a small ball $B$ and recall that $D \subseteq \Omega$ is such that $D+B \subseteq \Omega$, so that
$\phi(\cdot - x)$ lives on $\Omega$ whenever $x \in D$. We now have the tiling condition
\begin{equation} \label{tiling} \sum_{(a,b) \in S} {|V_g\phi((a,b)-(x,y))|}^2 = 1,\ \ \ \text{for $(x,y) \in D \times {\Bbb R}^d$}. \end{equation} 

In words, again,  the  $S$-time-frequency translations  of the function ${|V_g\phi|}^2$ is a tiling on $D \times {\Bbb R}^d$ (and a packing everywhere by \ref{packing}). 

Now $V_g\phi(t,\nu) = 0$ if $t \notin K-B$ so the tiling \eqref{tiling} must be effected
by the translates in 
$$S \cap \left((D+K+B) \times {\Bbb R}^d \right) .$$

One can take $\Omega$ big enough such that in addition   $K+D\subset \Omega$. This implies that  

$$S \cap \left((D+K+B) \times {\Bbb R}^d \right) \subseteq  S\cap (\Omega+B)\times {\Bbb R}^d .$$

Because of the (global) packing condition \eqref{packing} the size of the set
$S \cap \left((\Omega+B) \times (y+B_R)\right)$ is uniformly bounded in $y$ for any fixed  radius $R$.
So if
$$ \pi_2(S \cap \left((\Omega+B) \times {\Bbb R}^d) \right)$$
has zero lower density then so has the set $S \cap \left((\Omega+B) \times {\Bbb R}^d\right)$, which is incompatible with the tiling condition in (\ref{tiling}).

\vskip.125in 

\subsection{Extraction of an (essentially) linear triple} 
\label{subsectionextract} 

\vskip.125in 

We are going to prove that there exist $(a_i,b_i)$, $i=1,2,3$, such that the following conditions hold: 

\begin{itemize} 

\item i) $|a_i-a_j| \leq \frac{1}{100r}$ for some $r \ge 10^6$. 
\item ii) $r \leq |b_i-b_j| \leq 2r$ if $i \not=j$. 
\item iii) $b_1,b_2,b_3$ are in an $\frac{1}{r}$-neighborhood of a line. 

\end{itemize} 

\vskip.125in 

Let $S$ be the putative Gabor spectrum, $\pi_1(S)$ denote the projection of $S$ onto the first $d$ variables and $\pi_2(S)$ the projection of $S$ onto the last $d$ variables. Assume without loss of generality that $\vec{0} \in \pi_1(S)$. Let $Q$ be a cube in ${\Bbb R}^d$ of unit side-length, and let $A_Q$ denote the elements of $\pi_1(S)$ such that the intersection of $K+a$ and $Q$ has a non-empty interior. By Theorem \ref{localizedgabor}, $\pi_2(\pi_1^{-1}(A_Q))$ has positive lower density. Let $\Lambda_R$ denote the intersection of $\pi_2(\pi_1^{-1}(A_Q))$ with a ball of radius $R$, with $R$ suitably large. It follows that $\# \Lambda_R \approx R^d$, where here and throughout, $X \approx Y$ if there exists $C>0$ such that $C^{-1}Y \leq X \leq CY$. 

Cover $A_Q$ by a finitely overlapping family of $\approx r^{d}$ balls of radius $r^{-1}$, where $10^6 \leq r \leq R$. This gives us a decomposition $A_Q=\cup_j A_Q^j$. By the pigeon-hole principle, there exists $j_0$ such that $\# \pi_2(\pi_1^{-1}(A^{j_0}_Q)) \cap \Lambda_R \ge cr^{-d}R^d$. Let $H$ be a $(d-1)$-dimensional subspace of ${\Bbb R}^d$ and let $L$ denote its orthogonal subspace. We can divide the points in $\pi_2(\pi_1^{-1}(A^{j_0}_Q)) \cap \Lambda_R$ into $\approx {\left( \frac{R}{r} \right)}^{d-1}$ finitely overlapping classes where the $(d-1)$-coordinates determined by $H$ are within $\frac{1}{r}$ of one another. Applying the pigeon-hole principle once again and taking, say, $r=\log(R)$, we see that the triplet $(a_i,b_i)$, $i=1,2,3$, satisfying the conditions above exists. 

\vskip.125in

\subsection{Orthogonality and asymptotic expansions of the Fourier transform} 

\vskip.125in 

Let $K$ be a symmetric convex body with a smooth boundary and everywhere non-vanishing Gaussian curvature. 
Let $K$ be as above written in the form $ K=\{x: \rho(x) \leq 1 \},$ where $\rho:=\rho_K$ is the norm    that defines $K$, often called {\it the Minkowski functional (gauge)}. Define the dual functional by the relation
\begin{equation} \label{dualgauge} \rho^{*}(\xi)=\sup_{x \in \partial K} x \cdot \xi, \end{equation} where $\partial K$ denotes the boundary of $K$. We shall need the following result due to Herz (\cite{H62}; see also \cite{GS58}). 

\begin{lemma}\label{stationaryphase}   Let $K$ and $\rho_\ast$ be as above. Given $\omega \in S^{d-1}$, let $\kappa(\omega)$ denote the Gaussian curvature of $\partial K$ at the (unique) point where the unit normal is $\omega$. Then 
\begin{equation} \label{herzformula} \widehat{\chi}_K(\xi)=\kappa^{-\frac{1}{2}} \left( \frac{\xi}{|\xi|} \right) 
\sin \left( 2 \pi \left( \rho^{*}(\xi)-\frac{d-1}{8}\right) \right){|\xi|}^{-\frac{d+1}{2}}+{\mathcal D}_K(\xi), \end{equation} where 
$$ |{\mathcal D}_K(\xi)| \leq C_K {|\xi|}^{-\frac{d+3}{2}}.$$ 
\end{lemma} 

\vskip.125in 

\subsection{Conclusion of the argument} Let $(a_i,b_i)$, $i=1,2,3$ be the triple of elements of $S$ extracted in subsection \ref{subsectionextract} above. Going back to (\ref{orthoformula}), we have 
\begin{equation} \label{orthoours} \int \chi_K(x-a_i) \chi_K(x-a_j) e^{2 \pi i x \cdot (b_j-b_i)} dx=0. \end{equation} 

Note that the left hand side of the formula (\ref{orthoours}) is the Fourier transform of the indicator function of the symmetric convex body 
$(K+a_i) \cap (K+a_j)$ evaluated at $b_j-b_i$. The boundary of this body consists of a piece $\partial K$ and its translate. They intersect transversely at a smooth $(d-2)$-dimensional surface. The boundary of $(K+a_i) \cap (K+a_j)$ is smooth away from this surface. 

It is clear that the $(d-1)$-dimensional subspace $H$ (and hence its orthogonal subspace $L$) in the construction in the subsection \ref{subsectionextract} can be chosen such that $b_i-b_j$ is normal to the boundary $(K+a_i) \cap (K+a_j)$ at a point where the boundary is smooth. Then by applying a partition of unity, we see that Lemma \ref{stationaryphase} applies and in combination with the orthogonality relation (\ref{orthoours}) and the construction in subsection \ref{subsectionextract} implies that 
\begin{equation} \label{almostzeroes} \left| \sin \left( 2 \pi \rho_{a_i-a_j}^{*}(b_i-b_j)-\pi \frac{d-1}{4} \right) \right| \leq Cr^{-1}. \end{equation} 

This implies that 
\begin{equation} \label{zeroprep} \left| \sin \left( 2 \pi \left( \rho^{*}(b_i-b_j)+(a_i-a_j) \cdot (b_i-b_j)-\frac{d-1}{8} \right) \right) \right| \leq Cr^{-1} \end{equation} by writing out $\rho_{a_i-a_j}^{*}$ using the properties of the Fourier transform and the structure of the boundary of $(K+a_i) \cap (K+a_j)$ described above. 

Recalling the location of the zeroes of the sine function and the construction in Subsection \ref{subsectionextract}, we see that (\ref{zeroprep}) implies that 
$$ \left| \pi_L(b_i)-\pi_L(b_j)(1+c_{ij})-\frac{k}{2}+\frac{d-1}{8} \right| \leq C'r^{-1},$$ where 
$$ |c_{ij}| \leq \frac{1}{100r}.$$ 

It follows that 
$$ |\pi_L(b_2)-\pi_L(b_1)|=\frac{k_{12}}{2}-\frac{d-1}{8}+e_{12},$$ 
$$  |\pi_L(b_3)-\pi_L(b_2)|=\frac{k_{23}}{2}-\frac{d-1}{8}+e_{23},$$ and 
$$  |\pi_L(b_3)-\pi_L(b_1)|=\frac{k_{13}}{2}-\frac{d-1}{8}+e_{13},$$ where 
$$ |e_{ij}| \leq \frac{1}{50}$$ if $r$ is sufficiently large. 

\vskip.125in 

This implies that 
$$ \frac{k_{13}}{2}+\frac{d-1}{8}+e_{13}=\frac{k_{12}+k_{23}}{2}+\frac{d-1}{8}+\frac{d-1}{8}+e_{12}+e_{23},$$ so 
$$ \frac{k_{13}-k_{12}-k_{23}}{2}=\frac{d-1}{8}+e_{12}+e_{23}-e_{23}.$$

We conclude that $\frac{d-1}{4}+2(e_{12}+e_{23}-e_{23})$ is an integer. Since $d \not=1 \mod 4$, $\frac{d-1}{4}$ is, at best, half an integer. By above, 
$$ |2(e_{12}+e_{23}-e_{23})| \leq \frac{3}{25}.$$ 

This leads to a contradiction and the proof is complete. 

\vskip.125in 

\begin{remark} Note that the contradiction disappears if $d=1\mod 4$ because the $\frac{d-1}{8}$ term in (\ref{herzformula}) takes the form $\frac{4m+1-1}{8}=\frac{m}{2}$ and thus gets absorbed into the $\frac{k}{2}$ terms that precedes it. \end{remark} 

\vskip.25in
\section{Open problems: Spectral sets and Gabor windows} 
\label{openproblems}

\vskip.125in 

We conclude this paper with a series of questions that arise naturally out of considerations that led to our main result. We have shown that if $K$ is a symmetric convex body in ${\Bbb R}^d$, $d \not=1 \mod 4$, $d>1$, where $\partial K$ is smooth and has non-vanishing Gaussian curvature, then the indicator function of $K$ cannot serve as the window function for an orthogonal Gabor basis on $L^2({\Bbb R}^d)$. It was shown in \cite{IKT01} that under these assumptions, $L^2(K)$ does not possess an orthogonal basis of exponentials. This leads us to the following question. 

\vskip.125in 

\begin{definition} We say that $E \subset {\Bbb R}^d$ is {\it spectral} if $L^2({\Bbb R}^d)$ possesses an orthogonal basis of exponentials, namely a basis of the form ${ \{e^{-2 \pi i x \cdot b} \}}_{b \in B}$. We shall refer to $B$ as a {\it spectrum} for $E$. \end{definition} 

\begin{question}\label{nonspectralgabor} Suppose that $E$ is a {\it non-spectral} subset of ${\Bbb R}^d$. Is it possible that $\mathcal G(\chi_E, S) = \{\chi_E(x-a) e^{-2 \pi i x \cdot b} \}_{(a,b) \in S}$ is an orthogonal basis for $L^2({\Bbb R}^d)$? 
\end{question}

\vskip.125in 

A natural example to investigate in the context of the Question \ref{nonspectralgabor} is that of a triangle. It is not spectral which follows, for instance, from a result due to Kolountzakis which proves that non-symmetric convex sets are never spectral (\cite{K99}). A similar question can be posed for sets that tile by translation but are not spectral (see \cite{KM06} and the references contained therein). 

\vskip.125in 

Conversely, it is natural to ask whether all spectral sets generate Gabor windows:

\vskip.125in 

\begin{question}\label{spectralgabor} Suppose that $E$ is a spectral subset of ${\Bbb R}^d$. Is it true that there exists $S \subset {\Bbb R}^{2d}$ such that $\mathcal G(\chi_E, S) = \{\chi_E(x-a) e^{-2 \pi i x \cdot b} \}_{(a,b) \in S}$ is an orthogonal basis for $L^2({\Bbb R}^d)$? 
\end{question}  

The answer to this question is automatically affirmative if $E$ tiles ${\Bbb R}^d$ by translation. In this case we take $S=A \times B$, where $A$ is a tiling set for $E$ and $B$ is a spectrum. However, examples due Tao (\cite{T04}), Kolountzakis, Matolcsi (\cite{KM06}, \cite{M05}) and others show that there exist spectral sets which do not tile by translation. It would be natural to investigate the Question \ref{spectralgabor} starting with those constructions. 

It is worth recalling in this context that a result due to Fuglede (\cite{Fug74}, see also \cite{BHM17}) shows that if $E \subset {\Bbb R}^d$ tiles by a lattice $L$, then $E$ is spectral with a spectrum given by the dual lattice $L^{*}$. Conversely, if $E$ is spectral with a spectrum $L$ which happens to be a lattice, then $E$ tiles by translation by the dual lattice $L^{*}$. It follows that if $E$ tiles by a lattice $L$, then $\mathcal G(\chi_E, L \times L^{*})$ is an orthogonal basis for $L^2({\Bbb R}^d)$. Similarly, if $E$ is spectral with a lattice spectrum $L$, then $\mathcal G(\chi_E, L^{*} \times L)$ is an orthogonal basis for $L^2({\Bbb R}^d)$. 

\vskip.125in 

A possibly even more basic question arises from the discussion above. As we have seen, the case when $S=A \times B$ is somewhat special in the context of Gabor bases, so it makes sense to try to understand to what extent it is prevalent. 

\begin{question}\label{alwaysseparable} Does there exist a window function $g$ such that $\mathcal G(g, S) = \{g(x-a) e^{-2 \pi i x \cdot b} \}_{(a,b) \in S}$ is not an orthogonal basis for $L^2({\Bbb R}^d)$ for any $S$ of the form $A \times B$, but is an orthogonal basis for $L^2({\Bbb R}^d)$ for some $S$ not of that form? 
\end{question} 

In this general direction, Han and Wang (\cite{HW04}) proved that if $S=M {\Bbb Z}^{2d}$, where $M$ is a matrix with rational entries and $det(M)=1$, then there exists a compactly supported $g \in L^2({\Bbb R}^d)$ such that $\mathcal G(g, S) = \{g(x-a) e^{-2 \pi i x \cdot b} \}_{(a,b) \in S}$ is an orthogonal basis for $L^2({\Bbb R}^d)$. However, it is not immediately clear whether $g \in L^2({\Bbb R}^d)$ can be chosen in this case in such a way that $\mathcal G(g, S)$ is not an orthogonal basis for $L^2({\Bbb R}^d)$ for any $S$ of the form $A \times B$. 
 
\vskip.25in 

\section*{Acknowledgments} 
The authors wish to thank Akram Aldroubi, Carlos Cabrelli, Rachel Greenfeld, Mihalis Kolountzakis, Ursula Molter, Nir Lev and Shahaf Nitzan for helpful discussions and suggestions. 

\bibliographystyle{amsplain}

\begin{thebibliography}{99}

\bibitem{AAK17} E. Agora, J. Antezana, and M. Kolountzakis, {\it Tiling functions and Gabor orthonormal basis}, (preprint), (https://arxiv.org/pdf/1704.02831.pdf), (2017). 

\bibitem{B88} G. Battle, {\it Heisenberg proof of the Balian-Low theorem}, Letters in Mathematical Physics \textbf{15}, (1988), 175-177. 

\bibitem{BHM17} D. Barbieri, E. Hernandez and A. Mayeli, {\it Lattice sub-tilings and frames in LCA groups}, C. R. Math. Acad. Sci. Paris \textbf{355} (2017), no. 2, 193-199.

\bibitem{CLMP17} C. Cabrelli, D. Lee, U. Molter and G. Pfander, {\it Time-frequency shift invariance of Gabor spaces generated by integer lattices}, (arXiv:1705.02495) (2017). 

\bibitem{D90} I. Daubechies, {\it The wavelet transform, time-frequency localization and signal analysis}, IEEE Transactions on Information Theory, volume \textbf{36}, no. 5, (1990). 

\bibitem{Fug74} B. Fuglede, {\it Commuting self-adjoint partial differential operators and a group theoretic problem}, J. Funct. Anal. \textbf{16} (1974), 101-121.

\bibitem{G01} K. Gr\"ochenig, {\it Foundations of time-frequency analysis}, Applied and Numerical Harmonic Analysis. Birkh\"user Boston, Inc., Boston, MA, (2001). 

\bibitem{GM13} K. Gr\"ochenig and E. Malinnikova, {\it Phase space localization of Riesz bases in $L^2({\Bbb R}^d)$}, Rev. Mat. Iberoam. \textbf{29} (2013), no. 1, 115-134.

\bibitem{GS58} I. Gelfand and G. Shilov, {\it Generalized Functions}, Vol. 1, Academic Press, (1958).

\bibitem{Gr14} K. Gr\"ochenig, {\it The mystery of Gabor frames}, J Fourier Anal Appl (2014) 20:865-895. 

\bibitem{HW04} D. Han, Y. Wang, {\it The existence of Gabor bases and frames} Contemporary Math., \textbf{345}, (2004), special issue, p. 183-192.

\bibitem{H62} C. Herz, {\it Fourier transforms related to convex sets}, Ann. of Math. (2) \textbf{75} (1962) 81-92. 

\bibitem{IKMN18} {A. Iosevich, M. Kolountzakis, A. Mayeli and S. Nitzan}, {\it Geometric measure theory and the existence of orthogonal Gabor bases}, (in preparation), (2018). 


\bibitem{IKP01} {A. Iosevich, N.~H. Katz and S. Pedersen}, {\it Fourier bases and a distance problem of Erd\"os}. \newblock {\em Math. Res. Lett.}, \textbf{6}, 251-255, (2001).

\bibitem{IKT01} {A. Iosevich, N.~H. Katz and T. Tao,} {\it Convex bodies with a point of curvature do not have Fourier bases}, Amer. J. Math., \textbf{123}, 115-120, (2001).

\bibitem{IP00} A. Iosevich and S. Pedersen, {\it How large are the spectral gaps?} Pacific J. Math. \textbf{192} (2000), no. 2, 307-314. 

\bibitem{IR03} A. Iosevich and M. Rudnev, {\it A combinatorial approach to orthogonal exponentials}, Int. Math. Res. Not. (2003), no. 50, 2671-2685. 




\bibitem{K18} M. Kolountzakis, {\it Personal communication}, (2018). 

\bibitem{KL96} M. Kolountzakis and J. Lagarias, {\it Structure of tilings of the line by a function}, Duke Math. J. \textbf{82} (1996), no. 3, 653-678.

\bibitem{KM06} M. Kolountzakis and M. Matolcsi, {\it Tiles with no spectra}, Forum Math. \textbf{18} (2006), no. 3, 519-528.

\bibitem{K99} M. Kolountzakis, {\it Non-symmetric convex domains have no basis of exponentials}, Illinois journal of mathematics 44(3),  March 1999

\bibitem{LW03} Y. Liu and Y. Wang, {\it The Uniformity of Non-Uniform Gabor Bases}, Advances in Computational
Mathematics February (2003), Volume 18, Issue 2, pp 345-355.

\bibitem{M05} M. Matolcsi, {\it Fuglede conjecture fails in dimension $4$}, Proc. Amer. Math. Soc., \textbf{133} (2005), no. 10, 3021-3026.

\bibitem{M80} P. McMullen, {\it Convex bodies which tile space by translation}, Mathematika \textbf{27} (1980), no. 1, 113-121.

\bibitem{NO17} S. Nitzan and J-F. Olsen, {\it Balian-Low type theorem in finite dimensions}, (preprint), (2017). 

\bibitem{Ramanathan-Steger94} J. Ramanathan, T. Steger, {\it Incompleteness of Sparse Coherent States},  Applied and Computational Harmonic Analysis, 148-153 (1995).

\bibitem{T04} T. Tao, {\it Fuglede's conjecture is false in 5 and higher dimensions}, Math. Res. Lett. \textbf{11} (2004), no. 2-3, 251-258.

\end{thebibliography}


\begin{dajauthors}
\begin{authorinfo}[iosevich]
  Alex Iosevich\\
  Department of Mathematics\\
  University of Rochester\\
  Rochester, NY\\
  iosevich@math.rochester.edu\\
\end{authorinfo}
\begin{authorinfo}[mayeli]
  Azita Mayeli\\
  Department of Mathematics and Computer Science\\
  City University of New York (CUNY)\\
  The Graduate Centre and Queensborough, NY\\
  amayeli@gc.cuny.edu
\end{authorinfo}
\end{dajauthors}

\end{document}